\documentclass[leqno]{amsart}
\usepackage[latin9]{inputenc}
\usepackage{amsthm}
\usepackage{amssymb}
\usepackage{esint}
\usepackage{hyperref}

\makeatletter
\numberwithin{equation}{section}
\numberwithin{figure}{section}
\theoremstyle{plain}
\newtheorem{thm}{\protect\theoremname}
  \theoremstyle{plain}
  \newtheorem{lem}[thm]{\protect\lemmaname}
  \theoremstyle{remark}
  \newtheorem*{rem*}{\protect\remarkname}

\usepackage[numbers,sort&compress,square]{natbib}

\makeatother

  \providecommand{\lemmaname}{Lemma}
  \providecommand{\remarkname}{Remark}
\providecommand{\theoremname}{Theorem}

\begin{document}

\title[Identities for Carlitz $q$-Bernoulli polynomials]{Some identities of $q$-Bernoulli polynomials under symmetry group $S_3$}

\author{Dmitry V. Dolgy}
\address{Institute of Mathematics and Computer Sciences, Far Eastern Federal University, Vladivostok 690060, Russia}
\email{$d_{-}dol@mail.ru$}

\author{Dae San Kim}
\address{Department of Mathematics, Sogang University, Seoul 121-742, Republic
of Korea}
\email{dskim@sogang.ac.kr}

\author{Taekyun Kim}

\address{Department of Mathematics, Kwangwoon University, Seoul 139-701, Republic
of Korea}
\email{tkkim@kw.ac.kr}


\begin{abstract}
In this paper, we give some new identities of Carlitz $q$-Bernoulli
polynomials under symmetry group $S_3$.
The derivatives of identities are based on the $q$-Volkenborn integral
expression of the generating function for the Carlitz $q$-Bernoulli
polynomials and the $q$-Volkenborn integral equations that can be
expressed as the exponential generating functions for the $q$-power
sums.
\end{abstract}
\maketitle
\global\long\def\Zp{\mathbb{Z}_{p}}

\section{Introduction}

Let $p$ be a fixed prime number. Throughout this paper, $\Zp$, $\mathbb{Q}_{p}$
and $\mathbb{C}_{p}$ will, respectively, denote the ring of $p$-adic
integers, the field of $p$-adic rational numbers and the completion
of algebraic closure of $\mathbb{Q}_{p}$. Let $\nu_{p}$ be the normalized
exponential valuation of $\mathbb{C}_{p}$ with $\left|p\right|_{p}=p^{-\nu_{p}\left(p\right)}=\frac{1}{p}$
and let $q$ be an indeterminate in $\mathbb{C}_{p}$ with $\left|1-q\right|_{p}<p^{-\frac{1}{p-1}}$.
The $q$-number of $x$ is defined as $\left[x\right]_{q}=\frac{1-q^{x}}{1-q}$.
Let $UD\left(\Zp\right)$ be the space of uniformly differentiable
functions on $\Zp$. For $f\in UD\left(\Zp\right)$, the $q$-Volkenborn
integral on $\Zp$ is defined by Kim to be

\begin{equation}
I_{q}\left(f\right)=\int_{\Zp}f\left(x\right)d\mu_{q}\left(x\right)=\lim_{N\rightarrow\infty}\frac{1}{\left[p^{N}\right]_{q}}\sum_{x=0}^{p^{N}-1}f\left(x\right)q^{x},\quad\left(\mbox{see \cite{key-13}}\right).\label{eq:1}
\end{equation}

In (\ref{eq:1}), we note that
\begin{equation}
qI_{q}\left(f_{1}\right)-I_{q}\left(f\right)=\frac{q-1}{\log q}f^{\prime}\left(0\right)+\left(q-1\right)f\left(0\right),\label{eq:2}
\end{equation}
where $f_{1}\left(x\right)=f\left(x+1\right)$.

In general, one derives
\begin{equation}
q^{n}I_{q}\left(f_{n}\right)-I_{q}\left(f\right)=\frac{q-1}{\log q}\sum_{l=0}^{n-1}f^{\prime}\left(l\right)q^{l}+\left(q-1\right)\sum_{l=0}^{n-1}f\left(l\right)q^{l},\label{eq:3}
\end{equation}
where $f_{n}\left(x\right)=f\left(x+n\right)$, $\left(n\ge0\right)$,
(see \cite{key-14,key-15,key-25}).

It is well known that the Bernoulli numbers are given by
\begin{equation}
B_{0}=1,\quad\left(B+1\right)^{n}-B_{n}=\delta_{1,n},\quad\mbox{(see \cite{key-1,key-2,key-3,key-4,key-5,key-6,key-7,key-8,key-9,key-10,key-11,key-12,key-13,key-14,key-15,key-16,key-17,key-18,key-19,key-20,key-21,key-22,key-23,key-24,key-25,key-26,key-27})}\label{eq:4}
\end{equation}
with the usual convention about replacing $B^{n}$ by $B_{n}$.

The Bernoulli polynomials are defined by
\begin{equation}
B_{n}\left(x\right)=\sum_{l=0}^{n}\dbinom{n}{l}B_{l}x^{n-l},\quad\left(n\ge0\right),\quad\left(\mbox{see \cite{key-29,key-30}}\right).\label{eq:5}
\end{equation}

In \cite{key-4,key-5}, L. Carlitz considered a $q$-analogue of Bernoulli
numbers as follows :
\begin{equation}
\beta_{0,q}=1,\quad q\left(q\beta_{q}+1\right)^{n}-\beta_{n,q}=\begin{cases}
1 & \mbox{if }n=1,\\
0 & \mbox{if }n>1,
\end{cases}\label{eq:6}
\end{equation}
with the usual convention about replacing $\beta_{q}^{i}$ by $\beta_{i,q}$.

He also defined $q$-Bernoulli polynomials as follows :
\begin{equation}
\beta_{n,q}\left(x\right)=\sum_{l=0}^{n}\dbinom{n}{l}q^{lx}\beta_{l,q}\left[x\right]_{q}^{n-l},\quad\left(\mbox{see \cite{key-4,key-5}}\right).\label{eq:7}
\end{equation}

From (\ref{eq:6}), Carlitz derived the following equation :
\begin{equation}
\beta_{n,q}=\frac{1}{\left(1-q\right)^{n}}\sum_{l=0}^{n}\dbinom{n}{l}\left(-1\right)^{l}\frac{l+1}{\left[l+1\right]_{q}},\label{eq:8}
\end{equation}
and
\begin{equation}
\beta_{n,q}\left(x\right)=\frac{1}{\left(1-q\right)^{n}}\sum_{l=0}^{n}\dbinom{n}{l}\left(-1\right)^{l}q^{lx}\frac{l+1}{\left[l+1\right]_{q}},\quad\left(\mbox{see \cite{key-5}}\right).\label{eq:9}
\end{equation}

Carlitz $q$-Bernoulli numbers and polynomials are also given by $q$-Volkenborn
integrals on $\Zp$ due to T. Kim \cite{key-13} :
\begin{eqnarray}
\sum_{n=0}^{\infty}\beta_{n,q}\frac{t^{n}}{n!} & = & \int_{\Zp}e^{\left[y\right]_{q}t}d\mu_{q}\left(y\right)\label{eq:10}\\
 & = & \sum_{m=0}^{\infty}q^{m}e^{\left[m\right]_{q}t}\left(1-q-q^{m}t\right),\nonumber
\end{eqnarray}
and
\begin{align}
\sum_{n=0}^{\infty}\beta_{n,q}\left(x\right)\frac{t^{n}}{n!} & =\int_{\Zp}e^{\left[x+y\right]_{q}t}d\mu_{q}\left(y\right)\label{eq:11}\\
 & =\sum_{m=0}^{\infty}q^{m}e^{\left[x+m\right]_{q}t}\left(1-q-q^{x+m}t\right).\nonumber
\end{align}

The purpose of this paper is to give some new identities of Carlitz
$q$-Bernoulli polynomials under symmetry group $S_3$. The derivations of identities are based on the
$q$-Volkenborn integral expression of the generating function for
the Carlitz $q$-Bernoulli polynomials and the $q$-Volkenborn integrals
equations that can be expressed as the exponential generating functions
for the $q$-power sums.

\section{Some identities of Carlitz $q$-Bernoulli polynomials}

$\,$

By (\ref{eq:7}), we easily get
\begin{align}
\beta_{n,q}\left(x+y\right) & =\sum_{l=0}^{n}\dbinom{n}{l}q^{lx}\beta_{l,q}\left(y\right)\left[x\right]_{q}^{n-l}\label{eq:12}\\
 & =\sum_{l=0}^{n}\dbinom{n}{l}q^{\left(n-l\right)x}\beta_{n-l,q}\left(y\right)\left[x\right]_{q}^{l}.\nonumber
\end{align}

On the other hand, Carlitz also introduced the expression of $q$-Bernoulli
polynomials $\beta_{n,q}^{\left(h,k\right)}\left(x\right)$ as follows
:

\begin{equation}
\beta_{n,q}^{\left(h,k\right)}\left(x\right)=\frac{1}{\left(1-q\right)^{n}}\sum_{j=0}^{n}\dbinom{n}{j}\left(-1\right)^{j}q^{jx}\frac{\left(j+h\right)_{k}}{\left[j+h\right]_{q,k}},\label{eq:13}
\end{equation}
where
\begin{equation}
\left[j+h\right]_{q,k}=\left[j+h\right]_{q}\left[j+h-1\right]_{q}\cdots\left[j+h-k+1\right]_{q},\label{eq:14}
\end{equation}
and
\begin{equation}
\left(j+h\right)_{k}=\left(j+h\right)\left(j+h-1\right)\cdots\left(j+h-k+1\right),\quad\left(\mbox{see \cite{key-4,key-13}}\right).\label{eq:15}
\end{equation}

T. Kim \cite{key-13} obtained the Witt-type formula for $\beta_{n,q}^{\left(h,k\right)}\left(x\right)$
as follows :

\begin{equation}
\beta_{n,q}^{\left(h,k\right)}\left(x\right)=\int_{\Zp}\cdots\int_{\Zp}q^{\sum_{l=1}^{k}\left(h-l\right)y_{l}}\left[x+y_{1}+\cdots+y_{k}\right]_{q}^{n}d\mu_{q}\left(y_{1}\right)\cdots d\mu_{q}\left(y_{k}\right).\label{eq:16}
\end{equation}

If $k=1$, then $\beta_{n,q}^{\left(h,1\right)}\left(x\right)$ will
be simply denoted by $\beta_{n,q}^{\left(h\right)}\left(x\right)$
so that
\begin{equation}
\beta_{n,q}^{\left(h\right)}\left(x\right)=\int_{\Zp}q^{\left(h-1\right)y}\left[x+y\right]_{q}^{n}d\mu_{q}\left(y\right).\label{eq:17}
\end{equation}

The following simple facts will be used over and over again :
\begin{equation}
\left[a+b\right]_{q}=\left[a\right]_{q}+q^{a}\left[b\right]_{q}.\label{eq:18}
\end{equation}

By (\ref{eq:18}), we easily see that
\begin{equation}
\left[a+b+c\right]_{q}=\left[a\right]_{q}+q^{a}\left[b\right]_{q}+q^{a+b}\left[c\right]_{q},\label{eq:19}
\end{equation}
and
\begin{equation}
\left[ab\right]_{q}=\left[a\right]_{q}\left[b\right]_{q^{a}}.\label{eq:20}
\end{equation}

First, we will consider the following triple integral which is obviously
invariant under any permutations of $w_{1}$, $w_{2}$, $w_{3}$.

So the expression obtained from this after integration will also be
invariant under any permutations of $w_{1}$, $w_{2}$, $w_{3}$.
This observation is simple enough but it is the philosophy that underlies
this paper.

\begin{align}
I= & \int_{\mathbb{Z}_{p}^{3}}e^{\left[w_{2}w_{3}x_{1}+w_{1}w_{3}x_{2}+w_{1}w_{2}x_{3}+w_{1}w_{2}w_{3}\left(y_{1}+y_{2}+y_{3}\right)\right]_{q}t}\label{eq:21}\\
 & \times d\mu_{q^{w_{2}w_{3}}}\left(x_{1}\right)d\mu_{q^{w_{1}w_{3}}}\left(x_{2}\right)d\mu_{q^{w_{1}w_{2}}}\left(x_{3}\right).\nonumber
\end{align}

It is easy to show that
\begin{align}
 & \left[w_{2}w_{3}x_{1}+w_{1}w_{3}x_{2}+w_{1}w_{2}x_{3}+w_{1}w_{2}w_{3}\left(y_{1}+y_{2}+y_{3}\right)\right]_{q}\label{eq:22}\\
= & \left[w_{2}w_{3}\right]_{q}\left[x_{1}+w_{1}y_{1}\right]_{q^{w_{2}w_{3}}}+q^{w_{2}w_{3}\left(x_{1}+w_{1}y_{1}\right)}\left[w_{1}w_{3}\right]_{q}\left[x_{2}+w_{2}y_{2}\right]_{q^{w_{1}w_{3}}}\nonumber \\
 & +q^{w_{2}w_{3}\left(x_{1}+w_{1}y_{1}\right)+w_{1}w_{3}\left(x_{2}+w_{2}y_{2}\right)}\left[w_{1}w_{2}\right]_{q}\left[x_{3}+w_{3}y_{3}\right]_{q^{w_{1}w_{2}}}.\nonumber
\end{align}

So the integrand is
\begin{align}
 & e^{\left[w_{2}w_{3}x_{1}+w_{1}w_{3}x_{2}+w_{1}w_{2}x_{3}+w_{1}w_{2}w_{3}\left(y_{1}+y_{2}+y_{3}\right)\right]_{q}t}\label{eq:23}\\
= & e^{\left[w_{2}w_{3}\right]_{q}\left[x_{1}+w_{1}y_{1}\right]_{q^{w_{2}w_{3}}}t}e^{q^{w_{2}w_{3}\left(x_{1}+w_{1}y_{1}\right)}\left[w_{1}w_{3}\right]_{q}\left[x_{2}+w_{2}y_{2}\right]_{q^{w_{1}w_{3}}}t}\nonumber \\
 & \times e^{q^{w_{2}w_{3}\left(x_{1}+w_{1}y_{1}\right)+w_{1}w_{3}\left(x_{2}+w_{2}y_{2}\right)}\left[w_{1}w_{2}\right]_{q}\left[x_{3}+w_{3}y_{3}\right]_{q^{w_{1}w_{2}}}t}\nonumber \\
= & \sum_{n=0}^{\infty}\sum_{k+l+m=n}\dbinom{n}{k,l,m}\left[w_{2}w_{3}\right]_{q}^{k}\left[w_{1}w_{3}\right]_{q}^{l}\left[w_{1}w_{2}\right]_{q}^{m}q^{w_{1}w_{2}w_{3}\left(l+m\right)y_{1}}\nonumber \\
 & \times q^{w_{1}w_{2}w_{3}my_{2}}q^{w_{2}w_{3}\left(l+m\right)x_{1}}\left[x_{1}+w_{1}y_{1}\right]_{q^{w_{2}w_{3}}}^{k}q^{w_{1}w_{3}mx_{2}}\nonumber \\
 & \times\left[x_{2}+w_{2}y_{2}\right]_{q^{w_{1}w_{3}}}^{l}\left[x_{3}+w_{3}y_{3}\right]_{q^{w_{1}w_{2}}}^{m}\frac{t^{n}}{n!}.\nonumber
\end{align}

Thus the integral in (\ref{eq:21}) is given by
\begin{align}
I= & \sum_{n=0}^{\infty}\left\{ \sum_{k+l+m=n}\dbinom{n}{k,l,m}\left[w_{2}w_{3}\right]_{q}^{k}\left[w_{1}w_{3}\right]_{q}^{l}\left[w_{1}w_{2}\right]_{q}^{m}q^{w_{1}w_{2}w_{3}\left(l+m\right)y_{1}}\right.\label{eq:24}\\
 & \times q^{w_{1}w_{2}w_{3}my_{2}}\int_{\Zp}q^{w_{2}w_{3}\left(l+m\right)x_{1}}\left[x_{1}+w_{1}y_{1}\right]_{q^{w_{2}w_{3}}}^{k}d\mu_{q^{w_{2}w_{3}}}\left(x_{1}\right)\nonumber \\
 & \times\int_{\Zp}q^{w_{1}w_{3}mx_{2}}\left[x_{2}+w_{2}y_{2}\right]_{q^{w_{1}w_{3}}}^{l}d\mu_{q^{w_{1}w_{3}}}\left(x_{2}\right)\nonumber \\
 & \left.\times\int_{\Zp}\left[x_{3}+w_{3}y_{3}\right]_{q^{w_{1}w_{2}}}^{m}d\mu_{q^{w_{1}w_{2}}}\left(x_{3}\right)\right\} \frac{t^{n}}{n!}\nonumber \\
= & \sum_{n=0}^{\infty}\left\{ \sum_{k+l+m=n}\dbinom{n}{k,l,m}\left[w_{2}w_{3}\right]_{q}^{k}\left[w_{1}w_{3}\right]_{q}^{l}\right.\nonumber \\
 & \times\left[w_{1}w_{2}\right]_{q}^{m}q^{w_{1}w_{2}w_{3}\left(l+m\right)y_{1}}q^{w_{1}w_{2}w_{3}my_{2}}\nonumber \\
 & \left.\times\beta_{k,q^{w_{2}w_{3}}}^{\left(l+m+1\right)}\left(w_{1}y_{1}\right)\beta_{l,q^{w_{1}w_{3}}}^{\left(m+1\right)}\left(w_{2}y_{2}\right)\beta_{m,q^{w_{1}w_{2}}}\left(w_{3}y_{3}\right)\right\} \frac{t^{n}}{n!}.\nonumber
\end{align}

Thus, by (\ref{eq:24}), we get the following theorem.
\begin{thm}
\label{thm:1} Let $w_{1},w_{2},w_{3}$ be any positive integers,
$n$ any nonnegative integer. Then the following expression is invariant
under any permutation of $w_{1}$, $w_{2}$, $w_{3}$ so that it gives
us six symmetries :

\begin{align*}
 & \sum_{k+l+m=n}\dbinom{n}{k,l,m}\left[w_{2}w_{3}\right]_{q}^{k}\left[w_{1}w_{3}\right]_{q}^{l}\left[w_{1}w_{2}\right]_{q}^{m}q^{w_{1}w_{2}w_{3}\left(l+m\right)y_{1}}\\
 & \times q^{w_{1}w_{2}w_{3}my_{2}}\beta_{k,q^{w_{2}w_{3}}}^{\left(l+m+1\right)}\left(w_{1}y_{1}\right)\beta_{l,q^{w_{1}w_{3}}}^{\left(m+1\right)}\left(w_{2}y_{2}\right)\beta_{m,q^{w_{1}w_{2}}}\left(w_{3}y_{3}\right)\\
= & \sum_{k+l+m=n}\dbinom{n}{k,l,m}\left[w_{1}w_{3}\right]_{q}^{k}\left[w_{2}w_{3}\right]_{q}^{l}\left[w_{1}w_{2}\right]_{q}^{m}q^{w_{1}w_{2}w_{3}\left(l+m\right)y_{1}}\\
 & \times q^{w_{1}w_{2}w_{3}my_{2}}\beta_{k,q^{w_{1}w_{3}}}^{\left(l+m+1\right)}\left(w_{2}y_{1}\right)\beta_{l,q^{w_{2}w_{3}}}^{\left(m+1\right)}\left(w_{1}y_{2}\right)\beta_{m,q^{w_{1}w_{2}}}\left(w_{3}y_{3}\right)\\
= & \sum_{k+l+m=n}\dbinom{n}{k,l,m}\left[w_{1}w_{3}\right]_{q}^{k}\left[w_{1}w_{2}\right]_{q}^{l}\left[w_{2}w_{3}\right]_{q}^{m}q^{w_{1}w_{2}w_{3}\left(l+m\right)y_{1}}\\
 & \times q^{w_{1}w_{2}w_{3}my_{2}}\beta_{k,q^{w_{1}w_{3}}}^{\left(l+m+1\right)}\left(w_{2}y_{1}\right)\beta_{l,q^{w_{1}w_{2}}}^{\left(m+1\right)}\left(w_{3}y_{2}\right)\beta_{m,q^{w_{2}w_{3}}}\left(w_{1}y_{3}\right)\\
= & \sum_{k+l+m=n}\dbinom{n}{k,l,m}\left[w_{2}w_{3}\right]_{q}^{k}\left[w_{1}w_{2}\right]_{q}^{l}\left[w_{1}w_{3}\right]_{q}^{m}q^{w_{1}w_{2}w_{3}\left(l+m\right)y_{1}}\\
 & \times q^{w_{1}w_{2}w_{3}my_{2}}\beta_{k,q^{w_{2}w_{3}}}^{\left(l+m+1\right)}\left(w_{1}y_{1}\right)\beta_{l,q^{w_{1}w_{2}}}^{\left(m+1\right)}\left(w_{3}y_{2}\right)\beta_{m,q^{w_{1}w_{3}}}\left(w_{2}y_{3}\right)\\
= & \sum_{k+l+m=n}\dbinom{n}{k,l,m}\left[w_{1}w_{2}\right]_{q}^{k}\left[w_{2}w_{3}\right]_{q}^{l}\left[w_{1}w_{3}\right]_{q}^{m}q^{w_{1}w_{2}w_{3}\left(l+m\right)y_{1}}\\
 & \times q^{w_{1}w_{2}w_{3}my_{2}}\beta_{k,q^{w_{1}w_{2}}}^{\left(l+m+1\right)}\left(w_{3}y_{1}\right)\beta_{l,q^{w_{2}w_{3}}}^{\left(m+1\right)}\left(w_{1}y_{2}\right)\beta_{m,q^{w_{1}w_{3}}}\left(w_{2}y_{3}\right)\\
= & \sum_{k+l+m=n}\dbinom{n}{k,l,m}\left[w_{1}w_{2}\right]_{q}^{k}\left[w_{1}w_{3}\right]_{q}^{l}\left[w_{2}w_{3}\right]_{q}^{m}q^{w_{1}w_{2}w_{3}\left(l+m\right)y_{1}}\\
 & \times q^{w_{1}w_{2}w_{3}my_{2}}\beta_{k,q^{w_{1}w_{2}}}^{\left(l+m+1\right)}\left(w_{3}y_{1}\right)\beta_{l,q^{w_{1}w_{3}}}^{\left(m+1\right)}\left(w_{2}y_{2}\right)\beta_{m,q^{w_{2}w_{3}}}\left(w_{1}y_{3}\right).
\end{align*}

\end{thm}
We define, for nonnegative integers $n,m,w,T_{n,m}\left(w|q\right)$
as
\begin{equation}
T_{n,m}\left(w|q\right)=\sum_{i=0}^{w}q^{ni}\left[i\right]_{q}^{m}.\label{eq:25}
\end{equation}

In particular, for $w=0$, we have
\begin{equation}
T_{n,m}\left(0|q\right)=\begin{cases}
1, & \mbox{if }m=0\\
0, & \mbox{if }m>0,
\end{cases}\label{eq:26}
\end{equation}
and
\begin{equation}
T_{n,0}\left(w|q\right)=\begin{cases}
w+1 & \mbox{if }n=0\\
\left[w+1\right]_{q^{n}} & \mbox{if }n>0.
\end{cases}\label{eq:27}
\end{equation}

From (\ref{eq:3}), we have
\begin{align}
 & \left(q^{w_{1}w_{2}}\right)^{w_{3}}\int_{\Zp}e^{\left[w_{1}w_{2}\left(x+w_{3}\right)\right]_{q}t}d\mu_{q^{w_{1}w_{2}}}\left(x\right)\label{eq:28}\\
 & -\int_{\Zp}e^{\left[w_{1}w_{2}x\right]_{q}t}d\mu_{q^{w_{1}w_{2}}}\left(x\right)\nonumber \\
= & t\left[w_{1}w_{2}\right]_{q}\sum_{i=0}^{w_{3}-1}q^{2w_{1}w_{2}i}e^{\left[w_{1}w_{2}i\right]_{q}t}\nonumber \\
 & +\left(q-1\right)\left[w_{1}w_{2}\right]_{q}\sum_{i=0}^{w_{3}-1}q^{w_{1}w_{2}i}e^{\left[w_{1}w_{2}i\right]_{q}t}\nonumber \\
= & \sum_{m=0}^{\infty}T_{2,m}\left(w_{3}-1|q^{w_{1}w_{2}}\right)\left[w_{1}w_{2}\right]_{q}^{m+1}\frac{t^{m+1}}{m!}\nonumber \\
 & +\left(q-1\right)\sum_{m=0}^{\infty}T_{1,m}\left(w_{3}-1|q^{w_{1}w_{2}}\right)\left[w_{1}w_{2}\right]_{q}^{m+1}\frac{t^{m}}{m!}.\nonumber
\end{align}

Thus, we have the following lemma.
\begin{lem}
\label{lem:2} For $w_{1},w_{2},w_{3}\ge1$, we have
\begin{align*}
 & q^{w_{1}w_{2}w_{3}}\int_{\Zp}e^{\left[w_{1}w_{2}\left(x+w_{3}\right)\right]_{q}t}d\mu_{q^{w_{1}w_{2}}}\left(x\right)-\int_{\Zp}e^{\left[w_{1}w_{2}x\right]_{q}t}d\mu_{q^{w_{1}w_{2}}}\left(x\right)\\
= & \sum_{m=0}^{\infty}\frac{\left[w_{1}w_{2}\right]_{q}^{m}}{m!}t^{m}\int_{\Zp}\left(q^{w_{1}w_{2}w_{3}}\left[x+w_{3}\right]_{q^{w_{1}w_{2}}}^{m}-\left[x\right]_{q^{w_{1}w_{2}}}^{m}\right)d\mu_{q^{w_{1}w_{2}}}\left(x\right)\\
= & \sum_{m=0}^{\infty}T_{2,m}\left(w_{3}-1|q^{w_{1}w_{2}}\right)\left[w_{1}w_{2}\right]_{q}^{m+1}\frac{t^{m+1}}{m!}\\
 & +\left(q-1\right)\sum_{m=0}^{\infty}T_{1,m}\left(w_{3}-1|q^{w_{1}w_{2}}\right)\left[w_{1}w_{2}\right]_{q}^{m+1}\frac{t^{m}}{m!}\\
= & t\left[w_{1}w_{2}\right]_{q}\sum_{i=0}^{w_{3}-1}q^{2w_{1}w_{2}i}e^{\left[w_{1}w_{2}i\right]_{q}t}+\left(q-1\right)\left[w_{1}w_{2}\right]_{q}\sum_{i=0}^{w_{3}-1}q^{w_{1}w_{2}i}e^{\left[w_{1}w_{2}i\right]_{q}t}.
\end{align*}

\end{lem}
Now, we consider the following difference of triple integrals.

\begin{align}
I_{1}= & q^{w_{1}w_{2}w_{3}}\int_{\mathbb{Z}_{p}^{3}}e^{\left[w_{2}w_{3}x_{1}+w_{1}w_{3}x_{2}+w_{1}w_{2}x_{3}+w_{1}w_{2}w_{3}\left(y_{1}+y_{2}+1\right)\right]_{q}t}\label{eq:29}\\
 & \times d\mu_{q^{w_{2}w_{3}}}\left(x_{1}\right)d\mu_{q^{w_{1}w_{3}}}\left(x_{2}\right)d\mu_{q^{w_{1}w_{2}}}\left(x_{3}\right)\nonumber \\
 & -\int_{\mathbb{Z}_{p}^{3}}e^{\left[w_{2}w_{3}x_{1}+w_{1}w_{3}x_{2}+w_{1}w_{2}x_{3}+w_{1}w_{2}w_{3}\left(y_{1}+y_{2}\right)\right]_{q}t}\nonumber \\
 & \times d\mu_{q^{w_{2}w_{3}}}\left(x_{1}\right)d\mu_{q^{w_{1}w_{3}}}\left(x_{2}\right)d\mu_{q^{w_{1}w_{2}}}\left(x_{3}\right),\nonumber
\end{align}
which is obviously invariant under any permutations of $w_{1},w_{2},w_{3}$.

We put
\begin{equation}
a=a\left(x_{1}\right)=q^{w_{2}w_{3}\left(x_{1}+w_{1}y_{1}\right)},\quad b=b\left(x_{2}\right)=q^{w_{1}w_{3}\left(x_{2}+w_{2}y_{2}\right)}.\label{eq:30}
\end{equation}

Then, by (\ref{eq:29}), we get
\begin{align}
I_{1}= & \sum_{k,l=0}^{\infty}\left[w_{2}w_{3}\right]_{q}^{k}\left[w_{1}w_{3}\right]_{q}^{l}\frac{t^{k+l}}{k!l!}\int_{\mathbb{Z}_{p}^{2}}a^{l}\left[x_{1}+w_{1}y_{1}\right]_{q^{w_{2}w_{3}}}^{k}\left[x_{2}+w_{2}y_{2}\right]_{q^{w_{1}w_{3}}}^{l}\label{eq:31}\\
 & \times\left\{ \sum_{m=0}^{\infty}\frac{\left[w_{1}w_{2}\right]_{q}^{m}\left(abt\right)^{m}}{m!}\right.\nonumber \\
 & \times\left.\int_{\Zp}\left(q^{w_{1}w_{2}w_{3}}\left[x_{3}+w_{3}\right]_{q^{w_{1}w_{2}}}^{m}-\left[x_{3}\right]_{q^{w_{1}w_{2}}}^{m}\right)d\mu_{q^{w_{1}w_{2}}}\left(x_{3}\right)\right\} \nonumber \\
 & \times d\mu_{q^{w_{2}w_{3}}}\left(x_{1}\right)d\mu_{q^{w_{1}w_{3}}}\left(x_{2}\right).\nonumber
\end{align}

From Lemma \ref{lem:2}, the inner sum is

\begin{align}
 & \sum_{m=0}^{\infty}\frac{\left[w_{1}w_{2}\right]_{q}^{m+1}\left(abt\right)^{m+1}}{m!}T_{2,m}\left(w_{3}-1|q^{w_{1}w_{2}}\right)\label{eq:32}\\
+ & \left(q-1\right)\sum_{m=0}^{\infty}\frac{\left[w_{1}w_{2}\right]_{q}^{m+1}\left(abt\right)^{m}}{m!}T_{1,m}\left(w_{3}-1|q^{w_{1}w_{2}}\right).\nonumber
\end{align}

Thus, by (\ref{eq:31}) and (\ref{eq:32}), we get
\begin{align}
I_{1}= & \sum_{k,l,m=0}^{\infty}\left[w_{2}w_{3}\right]_{q}^{k}\left[w_{1}w_{3}\right]_{q}^{l}\left[w_{1}w_{2}\right]_{q}^{m+1}\frac{t^{k+l+m+1}}{k!l!m!}T_{2,m}\left(w_{3}-1|q^{w_{1}w_{2}}\right)\label{eq:33}\\
 & \times\int_{\Zp}a^{l+m+1}\left[x_{1}+m_{1}y_{1}\right]_{q^{w_{2}w_{3}}}^{k}d\mu_{q^{w_{2}w_{3}}}\left(x_{1}\right)\nonumber \\
 & \times\int_{\Zp}b^{m+1}\left[x_{2}+w_{2}y_{2}\right]_{q^{w_{1}w_{3}}}^{l}d\mu_{q^{w_{1}w_{3}}}\left(x_{2}\right)\nonumber
\end{align}
\begin{align*}
 & +\left(q-1\right)\sum_{k,l,m=0}^{\infty}\left[w_{2}w_{3}\right]_{q}^{k}\left[w_{1}w_{3}\right]_{q}^{l}\left[w_{1}w_{2}\right]_{q}^{m+1}\frac{t^{k+l+m}}{k!l!m!}\\
 & \times T_{1,m}\left(w_{3}-1|q^{w_{1}w_{2}}\right)\int_{\Zp}a^{l+m}\left[x_{1}+w_{1}y_{1}\right]_{q^{w_{2}w_{3}}}^{k}d\mu_{q^{w_{2}w_{3}}}\left(x_{1}\right)\\
 & \times\int_{\Zp}b^{m}\left[x_{2}+w_{2}y_{2}\right]_{q^{w_{1}w_{3}}}^{l}d\mu_{q^{w_{1}w_{3}}}\left(x_{2}\right).
\end{align*}

Recovering $a=q^{w_{2}w_{3}\left(x_{1}+w_{1}y_{1}\right)}$ and $b=q^{w_{1}w_{3}\left(x_{2}+w_{2}y_{2}\right)}$,
(\ref{eq:29}) can be rewritten as :
\begin{align}
I_{1}= & \sum_{n=0}^{\infty}\left\{ \sum_{k+l+m=n}\dbinom{n}{k,l,m}\left[w_{2}w_{3}\right]_{q}^{k}\left[w_{1}w_{3}\right]_{q}^{l}\left[w_{1}w_{2}\right]_{q}^{m+1}\right.\label{eq:34}\\
 & \times T_{2,m}\left(w_{3}-1|q^{w_{1}w_{2}}\right)q^{w_{1}w_{2}w_{3}\left(l+m+1\right)y_{1}}q^{w_{1}w_{2}w_{3}\left(m+1\right)y_{2}}\nonumber \\
 & \times\int_{\Zp}q^{w_{2}w_{3}\left(l+m+1\right)x_{1}}\left[x_{1}+w_{1}y_{1}\right]_{q^{w_{2}w_{3}}}^{k}d\mu_{q^{w_{2}w_{3}}}\left(x_{1}\right)\nonumber
\end{align}
\begin{align*}
 & \left.\times\int_{\Zp}q^{w_{1}w_{3}\left(m+1\right)x_{2}}\left[x_{2}+w_{2}y_{2}\right]_{q^{w_{1}w_{3}}}^{l}d\mu_{q^{w_{1}w_{3}}}\left(x_{2}\right)\right\} \frac{t^{n+1}}{n!}\\
 & +\left(q-1\right)\sum_{n=0}^{\infty}\left\{ \sum_{k+l+m=n}\dbinom{n}{k,l,m}\right.\\
 & \times\left[w_{2}w_{3}\right]_{q}^{k}\left[w_{1}w_{3}\right]_{q}^{l}\left[w_{1}w_{2}\right]_{q}^{m+1}T_{1,m}\left(w_{3}-1|q^{w_{1}w_{2}}\right)
\end{align*}
\begin{align*}
 & \times q^{w_{1}w_{2}w_{3}\left(l+m\right)y_{1}}q^{w_{1}w_{2}w_{3}my_{2}}\\
 & \times\int_{\Zp}q^{w_{2}w_{3}\left(l+m\right)x_{1}}\left[x_{1}+w_{1}y_{1}\right]_{q^{w_{2}w_{3}}}^{k}d\mu_{q^{w_{2}w_{3}}}\left(x_{1}\right)\\
 & \times\left.\int_{\Zp}q^{w_{1}w_{3}mx_{2}}\left[x_{2}+w_{2}y_{2}\right]_{q^{w_{1}w_{3}}}^{l}d\mu_{q^{w_{1}w_{3}}}\left(x_{2}\right)\right\} \frac{t^{n}}{n!}\\
= & \sum_{n=0}^{\infty}\left\{ \sum_{k+l+m=n}\dbinom{n}{k,l,m}\beta_{k,q^{w_{2}w_{3}}}^{\left(l+m+2\right)}\left(w_{1}y_{1}\right)\right.\\
 & \times\beta_{l,q^{w_{1}w_{3}}}^{\left(m+2\right)}\left(w_{2}y_{2}\right)T_{2,m}\left(w_{3}-1|q^{w_{1}w_{2}}\right)q^{w_{1}w_{2}w_{3}\left(l+m+1\right)y_{1}}\\
 & \left.\times q^{w_{1}w_{2}w_{3}\left(m+1\right)y_{2}}\left[w_{2}w_{3}\right]_{q}^{k}\left[w_{1}w_{3}\right]_{q}^{l}\left[w_{1}w_{2}\right]_{q}^{m+1}\right\} \frac{t^{n+1}}{n!}\\
 & +\left(q-1\right)\sum_{n=0}^{\infty}\left\{ \sum_{k+l+m=n}\dbinom{n}{k,l,m}\beta_{k,q^{w_{2}w_{3}}}^{\left(l+m+1\right)}\left(w_{1}y_{1}\right)\right.\\
 & \times\beta_{l,q^{w_{1}w_{3}}}^{\left(m+1\right)}\left(w_{2}y_{2}\right)T_{1,m}\left(w_{3}-1|q^{w_{1}w_{2}}\right)\\
 & \left.\times q^{w_{1}w_{2}w_{3}\left(l+m\right)y_{1}}q^{w_{1}w_{2}w_{3}my_{2}}\left[w_{2}w_{3}\right]_{q}^{k}\left[w_{1}w_{3}\right]_{q}^{l}\left[w_{1}w_{2}\right]_{q}^{m+1}\right\} \frac{t^{n}}{n!}.
\end{align*}

Separating the term corresponding to $n=0$ from the second term and
after rearranging, we get :
\begin{align}
I_{1}= & \left(q-1\right)\left[w_{1}w_{2}w_{3}\right]_{q}+\sum_{n=1}^{\infty}\left\{ \sum_{k+l+m=n-1}\dbinom{n}{k,l,m}\right.\label{eq:35}\\
 & \times\beta_{k,q^{w_{2}w_{3}}}^{\left(l+m+2\right)}\left(w_{1}y_{1}\right)\beta_{l,q^{w_{1}w_{3}}}^{\left(m+2\right)}\left(w_{2}y_{2}\right)T_{2,m}\left(w_{3}-1|q^{w_{1}w_{2}}\right)\nonumber \\
 & \times q^{w_{1}w_{2}w_{3}\left(l+m+1\right)y_{1}}q^{w_{1}w_{2}w_{3}\left(m+1\right)y_{2}}\left[w_{2}w_{3}\right]_{q}^{k}\left[w_{1}w_{3}\right]_{q}^{l}\left[w_{1}w_{2}\right]_{q}^{m+1}\nonumber
\end{align}
\begin{align*}
 & +\left(q-1\right)\sum_{k+l+m=n}\dbinom{n}{k,l,m}\beta_{k,q^{w_{2}w_{3}}}^{\left(l+m+1\right)}\left(w_{1}y_{1}\right)\\
 & \times\beta_{l,q^{w_{1}w_{3}}}^{\left(m+1\right)}\left(w_{2}y_{2}\right)T_{1,m}\left(w_{3}-1|q^{w_{1}w_{2}}\right)\\
 & \times q^{w_{1}w_{2}w_{3}\left(l+m\right)y_{1}}q^{w_{1}w_{2}w_{3}my_{2}}\\
 & \left.\times\left[w_{2}w_{3}\right]_{q}^{k}\left[w_{1}w_{3}\right]_{q}^{l}\left[w_{1}w_{2}\right]_{q}^{m+1}\right\} \frac{t^{n}}{n!}.
\end{align*}

As this expression is invariant under any permutations in $w_{1},w_{2},w_{3}$,
we get the following theorem.
\begin{thm}
\label{thm:3} Let $w_{1},w_{2},w_{3}\in\mathbb{Z}$ with $w_{1}\ge1$,
$w_{2}\ge1$, $w_{3}\ge1$. Then, for any positive integer $n$, the
following expressions
\begin{align*}
 & \sum_{k+l+m=n-1}\dbinom{n}{k,l,m}\beta_{k,q^{w_{\sigma\left(2\right)}w_{\sigma\left(3\right)}}}^{\left(l+m+2\right)}\left(w_{\sigma\left(1\right)}y_{1}\right)\\
 & \times\beta_{l,q^{w_{\sigma\left(1\right)}w_{\sigma\left(3\right)}}}^{\left(m+2\right)}\left(w_{\sigma\left(2\right)}y_{2}\right)T_{2,m}\left(w_{\sigma\left(3\right)}-1|q^{w_{\sigma\left(1\right)}w_{\sigma\left(2\right)}}\right)\\
 & \times q^{w_{1}w_{2}w_{3}\left(l+m+1\right)y_{1}}q^{w_{1}w_{2}w_{3}\left(m+1\right)y_{2}}\\
 & \times\left[w_{\sigma\left(2\right)}w_{\sigma\left(3\right)}\right]_{q}^{k}\left[w_{\sigma\left(1\right)}w_{\sigma\left(3\right)}\right]_{q}^{l}\left[w_{\sigma\left(1\right)}w_{\sigma\left(2\right)}\right]_{q}^{m+1}\\
 & +\left(q-1\right)\sum_{k+l+m=n}\dbinom{n}{k,l,m}\beta_{k,q^{w_{\sigma\left(2\right)}w_{\sigma\left(3\right)}}}^{\left(l+m+1\right)}\left(w_{\sigma\left(1\right)}y_{1}\right)\\
 & \times\beta_{l,q^{w_{\sigma\left(1\right)}w_{\sigma\left(3\right)}}}^{\left(m+1\right)}\left(w_{\sigma\left(2\right)}y_{2}\right)T_{1,m}\left(w_{\sigma\left(3\right)}-1|q^{w_{\sigma\left(1\right)}w_{\sigma\left(2\right)}}\right)\\
 & \times q^{w_{1}w_{2}w_{3}\left(l+m\right)y_{1}}q^{w_{1}w_{2}w_{3}my_{2}}\left[w_{\sigma\left(2\right)}w_{\sigma\left(3\right)}\right]_{q}^{k}\left[w_{\sigma\left(1\right)}w_{\sigma\left(3\right)}\right]_{q}^{l}\left[w_{\sigma\left(1\right)}w_{\sigma\left(2\right)}\right]_{q}^{m+1}
\end{align*}
are the same for any $\sigma\in S_{3}$. \end{thm}
\begin{rem*}
We can get interesting identities by letting $w_{3}=1$ or by letting
$w_{2}=w_{3}=1$, in view of (\ref{eq:26}). However, writing those
down requires much space. So we omit it.
\end{rem*}
With the same $a=q^{w_{2}w_{3}\left(x_{1}+w_{1}y_{1}\right)}$, $b=q^{w_{1}w_{3}\left(x_{2}+w_{2}y_{2}\right)}$
as in (\ref{eq:29}), $I_{1}$ can be rewritten as
\begin{align}
I_{1}= & \int_{\mathbb{Z}_{p}^{2}}e^{\left[w_{2}w_{3}\right]_{q}\left[x_{1}+w_{1}y_{1}\right]_{q^{w_{2}w_{3}}}t}e^{\left[w_{1}w_{3}\right]_{q}\left[x_{2}+w_{2}y_{2}\right]_{q^{w_{1}w_{3}}}\left(at\right)}\label{eq:36}\\
 & \times\left\{ \int_{\Zp}\left(q^{w_{1}w_{2}w_{3}}e^{\left[w_{1}w_{2}\left(x_{3}+w_{3}\right)\right]_{q}\left(abt\right)}-e^{\left[w_{1}w_{2}x_{3}\right]_{q}\left(abt\right)}\right)d\mu_{q^{w_{1}w_{2}}}\left(x_{3}\right)\right\} \nonumber \\
 & \times d\mu_{q^{w_{2}w_{3}}}\left(x_{1}\right)d\mu_{q^{w_{1}w_{3}}}\left(x_{2}\right).\nonumber
\end{align}

From Lemma \ref{lem:2} and (\ref{eq:36}), we note that the inner
integral is equal to
\begin{align}
 & abt\left[w_{1}w_{2}\right]_{q}\sum_{i=0}^{w_{3}-1}q^{2w_{1}w_{2}i}e^{\left[w_{1}w_{2}i\right]_{q}abt}\label{eq:37}\\
 & +\left(q-1\right)\left[w_{1}w_{2}\right]_{q}\sum_{i=0}^{w_{3}-1}q^{w_{1}w_{2}i}e^{\left[w_{1}w_{2}i\right]_{q}abt}.\nonumber
\end{align}

Thus, by (\ref{eq:36}) and (\ref{eq:37}), we get
\begin{align}
I_{1}= & t\left[w_{1}w_{2}\right]_{q}\sum_{i=0}^{w_{3}-1}q^{2w_{1}w_{2}i}\int_{\mathbb{Z}_{p}^{2}}abe^{\left[w_{2}w_{3}\right]_{q}\left[x_{1}+w_{1}y_{1}\right]_{q^{w_{2}w_{3}}}t}\label{eq:38}\\
 & \times e^{\left[w_{1}w_{3}\right]_{q}\left[x_{2}+w_{2}y_{2}+\frac{w_{2}}{w_{3}}i\right]_{q^{w_{1}w_{3}}}at}d\mu_{q^{w_{2}w_{3}}}\left(x_{1}\right)d\mu_{q^{w_{1}w_{3}}}\left(x_{2}\right)\nonumber \\
 & +\left(q-1\right)\left[w_{1}w_{2}\right]_{q}\sum_{i=0}^{w_{3}-1}q^{w_{1}w_{2}i}\int_{\mathbb{Z}_{p}^{2}}e^{\left[w_{2}w_{3}\right]_{q}\left[x_{1}+w_{1}y_{1}\right]_{q^{w_{2}w_{3}}}t}\nonumber \\
 & \times e^{\left[w_{1}w_{3}\right]_{q}\left[x_{2}+w_{2}y_{2}+\frac{w_{2}}{w_{3}}i\right]_{q^{w_{1}w_{3}}}at}d\mu_{q^{w_{2}w_{3}}}\left(x_{1}\right)d\mu_{q^{w_{1}w_{3}}}\left(x_{2}\right)\nonumber
\end{align}
\begin{align*}
= & t\left[w_{1}w_{2}\right]_{q}\sum_{i=0}^{w_{3}-1}q^{2w_{1}w_{2}i}\sum_{k,l=0}^{\infty}\frac{t^{k+l}}{k!l!}q^{w_{1}w_{2}w_{3}\left(l+1\right)y_{1}}q^{w_{1}w_{2}w_{3}y_{2}}\left[w_{2}w_{3}\right]_{q}^{k}\\
 & \times\left[w_{1}w_{3}\right]_{q}^{l}\int_{\Zp}q^{w_{2}w_{3}\left(l+1\right)x_{1}}\left[x_{1}+w_{1}y_{1}\right]_{q^{w_{2}w_{3}}}^{k}d\mu_{q^{w_{2}w_{3}}}\left(x_{1}\right)\\
 & \times\int_{\Zp}q^{w_{1}w_{3}x_{2}}\left[x_{2}+w_{2}y_{2}+\frac{w_{2}}{w_{3}}i\right]_{q^{w_{1}w_{3}}}^{l}d\mu_{q^{w_{1}w_{3}}}\left(x_{2}\right)
\end{align*}
\begin{align*}
 & +\left(q-1\right)\left[w_{1}w_{2}\right]_{q}\sum_{i=0}^{w_{3}-1}q^{w_{1}w_{2}i}\sum_{k,l=0}^{\infty}\frac{t^{k+l}}{k!l!}q^{w_{1}w_{2}w_{3}ly_{1}}\left[w_{2}w_{3}\right]_{q}^{k}\\
 & \times\left[w_{1}w_{3}\right]_{q}^{l}\int_{\Zp}q^{w_{2}w_{3}lx_{1}}\left[x_{1}+w_{1}y_{1}\right]_{q^{w_{2}w_{3}}}^{k}d\mu_{q^{w_{2}w_{3}}}\left(x_{1}\right)\\
 & \times\int_{\Zp}\left[x_{2}+w_{2}y_{2}+\frac{w_{2}}{w_{3}}i\right]_{q^{w_{1}w_{3}}}^{l}d\mu_{q^{w_{1}w_{3}}}\left(x_{2}\right)
\end{align*}
\begin{align*}
= & \sum_{n=0}^{\infty}\left\{ \sum_{k=0}^{n}\dbinom{n}{k}\beta_{k,q^{w_{2}w_{3}}}^{\left(n-k+2\right)}\left(w_{1}y_{1}\right)q^{w_{1}w_{2}w_{3}\left(n-k+1\right)y_{1}}\right.\\
 & \times q^{w_{1}w_{2}w_{3}y_{2}}\left[w_{_{1}}w_{2}\right]_{q}\left[w_{2}w_{3}\right]_{q}^{k}\\
 & \left.\times\left[w_{1}w_{3}\right]_{q}^{n-k}\sum_{i=0}^{w_{3}-1}q^{2w_{1}w_{2}i}\beta_{n-k,q^{w_{1}w_{3}}}^{\left(2\right)}\left(w_{2}y_{2}+\frac{w_{2}}{w_{3}}i\right)\right\} \frac{t^{n+1}}{n!}
\end{align*}
\begin{align*}
 & +\left(q-1\right)\sum_{n=0}^{\infty}\left\{ \sum_{k=0}^{n}\dbinom{n}{k}\beta_{k,q^{w_{2}w_{3}}}^{\left(n-k+1\right)}\left(w_{1}y_{1}\right)q^{w_{1}w_{2}w_{3}\left(n-k\right)y_{1}}\right.\\
 & \times\left[w_{1}w_{2}\right]_{q}\left[w_{2}w_{3}\right]_{q}^{k}\left[w_{1}w_{3}\right]_{q}^{n-k}\\
 & \left.\times\sum_{i=0}^{w_{3}-1}q^{w_{1}w_{2}i}\beta_{n-k,q^{w_{1}w_{3}}}\left(w_{2}y_{2}+\frac{w_{2}}{w_{3}}i\right)\right\} \frac{t^{n}}{n!}.
\end{align*}

Separating the term corresponding to $n=0$ from the second sum and
after rearranging, we obtain
\begin{align*}
I_{1}= & \left(q-1\right)\left[w_{1}w_{2}w_{3}\right]_{q}+\sum_{n=1}^{\infty}\left\{ n\sum_{k=0}^{n-1}\dbinom{n-1}{k}\beta_{k,q^{w_{2}w_{3}}}^{\left(n-k+1\right)}\left(w_{1}y_{1}\right)q^{w_{1}w_{2}w_{3}\left(n-k\right)y_{1}}\right.\\
 & \times q^{w_{1}w_{2}w_{3}y_{2}}\left[w_{1}w_{2}\right]_{q}\left[w_{2}w_{3}\right]_{q}^{k}\left[w_{1}w_{3}\right]_{q}^{n-1-k}\sum_{i=0}^{w_{3}-1}q^{2w_{1}w_{2}i}\\
 & \times\beta_{n-1-k,q^{w_{1}w_{3}}}^{\left(2\right)}\left(w_{2}y_{2}+\frac{w_{2}}{w_{3}}i\right)+\left(q-1\right)\sum_{k=0}^{n}\dbinom{n}{k}\beta_{k,q^{w_{2}w_{3}}}^{\left(n-k+1\right)}\left(w_{1}y_{1}\right)q^{w_{1}w_{2}w_{3}\left(n-k\right)y_{1}}\\
 & \left.\times\left[w_{1}w_{2}\right]_{q}\left[w_{2}w_{3}\right]_{q}^{k}\left[w_{1}w_{3}\right]_{q}^{n-k}\sum_{i=0}^{w_{3}-1}q^{w_{1}w_{2}i}\beta_{n-k,q^{w_{1}w_{3}}}\left(w_{2}y_{2}+\frac{w_{2}}{w_{3}}i\right)\right\} \frac{t^{n}}{n!}.
\end{align*}

As this expression is invariant under any permutations in $w_{1}$,
$w_{2}$, $w_{3}$, we have the following results.
\begin{thm}
\label{thm:4} Let $w_{1},w_{2},w_{3}\in\mathbb{Z}$ with $w_{1}\ge1$,
$w_{2}\ge1$, $w_{3}\ge1$. Then, for any positive integer $n$, the
following expressions
\begin{align*}
 & n\sum_{k=0}^{n-1}\dbinom{n-1}{k}\beta_{k,q^{w_{\sigma\left(2\right)}w_{\sigma\left(3\right)}}}^{\left(n-k+1\right)}\left(w_{\sigma\left(1\right)}y_{1}\right)q^{w_{1}w_{2}w_{3}\left(n-k\right)y_{1}}\\
 & \times q^{w_{1}w_{2}w_{3}y_{2}}\left[w_{\sigma\left(1\right)}w_{\sigma\left(2\right)}\right]_{q}\left[w_{\sigma\left(2\right)}w_{\sigma\left(3\right)}\right]_{q}^{k}\left[w_{\sigma\left(1\right)}w_{\sigma\left(3\right)}\right]_{q}^{n-1-k}\\
 & \times\sum_{i=0}^{w_{\sigma\left(3\right)}-1}q^{2w_{\sigma\left(1\right)}w_{\sigma\left(2\right)}i}\beta_{n-1-k,q^{w_{\sigma\left(1\right)}w_{\sigma\left(3\right)}}}^{\left(2\right)}\left(w_{\sigma\left(2\right)}y_{2}+\frac{w_{\sigma\left(2\right)}}{w_{\sigma\left(3\right)}}i\right)\\
 & +\left(q-1\right)\sum_{k=0}^{n}\dbinom{n}{k}\beta_{k,q^{w_{\sigma\left(2\right)}w_{\sigma\left(3\right)}}}^{\left(n-k+1\right)}\left(w_{\sigma\left(1\right)}y_{1}\right)q^{w_{1}w_{2}w_{3}\left(n-k\right)y_{1}}\\
 & \times\left[w_{\sigma\left(1\right)}w_{\sigma\left(2\right)}\right]_{q}\left[w_{\sigma\left(2\right)}w_{\sigma\left(3\right)}\right]_{q}^{k}\left[w_{\sigma\left(1\right)}w_{\sigma\left(3\right)}\right]_{q}^{n-k}\\
 & \times\sum_{i=0}^{w_{\sigma\left(3\right)}-1}q^{w_{\sigma\left(1\right)}w_{\sigma\left(2\right)}i}\beta_{n-k,q^{w_{\sigma\left(1\right)}w_{\sigma\left(3\right)}}}\left(w_{\sigma\left(2\right)}y_{2}+\frac{w_{\sigma\left(2\right)}}{w_{\sigma\left(3\right)}}i\right)
\end{align*}
are all the same for any $\sigma\in S_{3}$.\end{thm}
\begin{rem*}
Using (\ref{eq:25}) and by specializing $w_{3}=1$ or $w_{2}=w_{3}=1$,
we can obtain many interesting identities. However, as this requires
much space, we will omit those.
\end{rem*}

\section*{Acknowledgement}

This work was supported by the National Research Foundation of Korea(NRF)
grant funded by the Korea government(MOE) (No.2012R1A1A2003786 ) and partially supported by Kwangwoon University 2014.

\providecommand{\bysame}{\leavevmode\hbox to3em{\hrulefill}\thinspace}
\providecommand{\MR}{\relax\ifhmode\unskip\space\fi MR }
\providecommand{\MRhref}[2]{%
  \href{http://www.ams.org/mathscinet-getitem?mr=#1}{#2}
}
\providecommand{\href}[2]{#2}

\end{document}